\documentclass{amsart}
\usepackage{amssymb}
\usepackage{amsfonts}
\usepackage{latexsym}
\usepackage[all]{xy}
\usepackage{amscd}
\usepackage{comment}

\newtheorem{theorem}{Theorem}[section]
\newtheorem{lemma}[theorem]{Lemma}
\newtheorem{proposition}[theorem]{Proposition}
\newtheorem{corollary}[theorem]{Corollary}
\theoremstyle{definition}
\newtheorem{definition}[theorem]{Definition}
\newtheorem{example}[theorem]{Example}

\newtheorem{remark}[theorem]{Remark}

\newcommand{\C}{\mathcal{C}}
\newcommand{\D}{\mathcal{D}}
\newcommand{\E}{\mathcal{E}}

\newcommand{\M}{\mathcal{M}}

\newcommand{\R}{\mathcal{R}}

\renewcommand{\S}{\mathcal{S}}
\renewcommand{\O}{\mathcal{O}}
\newcommand{\id}{\text{id}}

\newcommand{\ot}{\otimes}

\newcommand{\Hom}{\mbox{Hom}}

\newcommand{\Rep}{\mbox{Rep}}

\newcommand{\Ker}{\mbox{Ker }}
\renewcommand{\Im}{\mbox{Im }}

\renewcommand{\Vec}{\mbox{Vec}}

\begin{document}

\title{Some properties of group-theoretical categories}

\author{Shlomo Gelaki}
\address{Department of Mathematics,
Technion - Israel Institute of Technology, Haifa 32000, Israel.}
\email{gelaki@math.technion.ac.il}

\author{Deepak Naidu}
\address{Department of Mathematics and Statistics,
University of New Hampshire, Durham, NH 03824, USA.}
\email{dnaidu@unh.edu}

\date{\today}

\begin{abstract}
We first show that every group-theoretical category is graded by a
certain double coset ring. As a consequence, we obtain a necessary
and sufficient condition for a group-theoretical category to be
nilpotent. We then give an explicit description of the simple objects in a
group-theoretical category (following \cite{O2}) and of the group
of invertible objects of a group-theoretical category, in
group-theoretical terms. Finally, under certain restrictive
conditions, we describe the universal grading group of a
group-theoretical category.
\end{abstract}

\maketitle

%%%%%%%%%%%%%%%%%%%%%%%%%%%%%%%%%%%%%%%%%%%%%%%%%%%%%%%%%%%%%%%%%%%%%
\begin{section}{Introduction}

Group-theoretical categories were introduced and studied in
\cite{ENO} and \cite{O1}. They constitute a fundamental class of
fusion categories which are defined, as the name suggests, by a
certain finite group data. For example, for a finite group $G$ its
representation category $\Rep(G)$ is group-theoretical. As an
indication of the centrality of group-theoretical categories in
the theory of fusion categories we mention the following
observation: all known complex semisimple Hopf algebras (as far as
we know) have group-theoretical representation categories. In
fact, it was asked in \cite{ENO} whether it is true that any
complex semisimple Hopf algebra is group-theoretical. It is thus
highly desirable to study group-theoretical categories and
understand as much as possible about them in the language of group
theory.

The notion of a nilpotent fusion category was introduced and
studied in \cite{GN}. For example, it is not hard to show that if
$G$ is a finite group then $\Rep(G)$ is nilpotent if and only if
$G$ is nilpotent. In \cite{DGNO} nilpotent modular categories are
studied, and in particular it is discussed when they are
group-theoretical. Therefore a very natural question arises: what
are necessary and sufficient conditions for a group-theoretical
category to be nilpotent? The answer to this question is one of the main
results of this paper (see Corollary \ref{main corollary}).

Other important invariants of a fusion category $\mathcal{C}$ are
its pointed subcategory $\mathcal{C}_{pt}$ (the subcategory
generated by the group of invertible objects in $\mathcal{C}$),
its adjoint subcategory $\mathcal{C}_{ad}$ \cite{ENO} and its
universal grading group $U(\mathcal{C})$ \cite{GN}. Descriptions
of $\mathcal{C}_{pt}$ for a general group-theoretical category
$\mathcal{C}$, and $\mathcal{C}_{ad}$, $U(\mathcal{C})$ for a
special class of group-theoretical categories are other results of
this paper (see
Theorem \ref{grinv} and Proposition
\ref{unigr}).

The organization of the paper is as follows. Section 2 contains
necessary preliminaries about fusion categories, module
categories, and group-theoretical categories. We also recall some
definitions from \cite{GN} concerning nilpotent fusion categories and
based rings.
We also recall some basic definitions and results from group
theory.

In Section 3 we introduce the notion of a fusion category graded by a
based ring. Let $H$ be a subgroup of a finite group $G$. We introduce
a based ring which we call {\em double coset ring} arising from the set
$H \backslash G/H$ of double cosets of $H$ in $G$. We give a necessary and
sufficient condition for the double coset ring to be nilpotent (see
Proposition \ref{nilpotent double coset}).

In Section 4 we first show that every group-theoretical category
is graded by a certain double coset ring. As a consequence, we
obtain a necessary and sufficient condition for a
group-theoretical category to be nilpotent.

In Section 5 we give an explicit description of the simple objects in a
group-theoretical category (following Proposition 3.2 in \cite{O2}; see Theorem
\ref{simpobj}) and of the group of invertible objects of a
group-theoretical category, in group-theoretical terms.

In Section 6, we describe the universal grading group of a
group-theoretical category, under certain restrictive conditions.

\end{section}
%%%%%%%%%%%%%%%%%%%%%%%%%%%%%%%%%%%%%%%
\par
{\bf Acknowledgments.} Part of this work was done while the first author was on
Sabbatical in the departments of mathematics at the University of
New Hampshire and MIT; he is grateful for their warm hospitality.
The research of the first author was partially supported by the
Israel Science Foundation (grant No. 125/05). The authors would
like to thank P. Etingof and D. Nikshych for useful discussions.

%%%%%%%%%%%%%%%%%%%%%%%%%%%%%%%%%%%%%%%%%%%%%%%%%%%%%%%%%%%%%%%%%%%%%%%
\begin{section}{Preliminaries}

%%%%%%%%%%%%%%%%%%%%%%%%%%%%%%%%%%%%%%%
\begin{subsection}{Fusion categories and their module categories}\hfill

Throughout this paper we work over an algebraically closed field
$k$ of characteristic $0$. All categories considered in this work
are assumed to be $k$-linear and semisimple with finite
dimensional $\Hom$-spaces and finitely many isomorphism classes of
simple objects. All functors are assumed to be additive and
$k$-linear. Unless otherwise stated all cocycles appearing in this
work will have coefficients in the trivial module $k^\times$.

A {\em fusion category} over $k$ is a
$k$-linear semisimple rigid tensor category with finitely many
isomorphism classes of simple objects and finite dimensional Hom-spaces
such that the neutral object is simple \cite{ENO}.

A fusion category is said to be {\em pointed} if all its simple
objects are invertible. A typical example of a pointed category
is $\Vec_G^{\omega}$ - the category of finite dimensional vector
spaces over $k$ graded by the finite group $G$. The morphisms in this category
are linear transformations that respect the grading and the associativity
constraint is given by the normalized $3$-cocycle $\omega$ on $G$.

Let $\C = (\C, \, \otimes, \, 1_{\C}, \, \alpha, \, \lambda, \, \rho)$ be a
tensor category, where $1_{\C}$, $\alpha$, $\lambda$, and $\rho$
are the unit object, the associativity constraint,
the left unit constraint, and the right unit constraint, respectively.
A right {\em module category} over $\C$ (see \cite{O1} and references therein) is
a category $\M$ together with an exact bifunctor $\otimes: \M \times \C \to
\M$ and natural isomorphisms
$\mu_{M, \, X, \, Y}: M \otimes (X \otimes Y) \to (M \otimes X) \otimes Y, \,\,
\tau_M: M \otimes 1_{\C} \to M$, for all $M \in \M, \, X, Y \in \C$, such that
the following two equations hold for all $M \in \M, \, X, Y, Z \in \C$:
\begin{equation*}
\label{module pentagon}
\mu_{M \otimes X, \, Y, \, Z} \, \circ
\, \mu_{M, \, X, \,  Y \otimes Z} \, \circ
\, (\id_M \otimes \alpha_{X,Y,Z})
= (\mu_{M, \, X, \, Y}\otimes \id_Z) \, \circ
\, \mu_{M, \, X \otimes Y, \, Z},
\end{equation*}
\begin{equation*}
\label{module triangle}
(\tau_M \otimes \id_Y) \, \circ
\, \mu_{M, \, 1_\C, \, Y}
= \id_M \otimes \lambda_Y.
\end{equation*}
Let $(\M_1, \, \mu^1, \tau^1)$ and $(\M_2, \, \mu^2, \tau^2)$ be two
right module categories over
$\C$. A $\C$-{\em module functor} from $\M_1$ to $\M_2$
is a functor $F: \M_1\to \M_2$ together with natural isomorphisms
$\gamma_{M, \, X}: F(M \otimes X) \to F(M) \otimes X$, for all
$M \in \M_1, \, X \in \C$, such that the following two
equations hold for all $M \in \M_1, \, X, Y \in \C$:
\begin{equation*}
\label{module functor pentagon}
(\gamma_{M, \, X} \otimes \id_Y) \, \circ
\, \gamma_{M\otimes X, \, Y} \, \circ
\, F(\mu^1_{M, \, X, \, Y})
= \mu^2_{F(M), \, X, \, Y} \, \circ
\, \gamma_{M, X \otimes Y},
\end{equation*}
\begin{equation*}
\label{module functor triangle}
\tau^1_{F(M)} \, \circ \, \gamma_{M, \, 1_\C} = F(\tau^1_M).
\end{equation*}
Two module categories $\M_1$ and $\M_2$ over $\C$ are {\em equivalent}
if there exists a module functor from $\M_1$ to $\M_2$ which is an
equivalence of categories.
For two module categories $\M_1$ and $\M_2$ over a tensor category
$\C$ their {\em direct sum} is the category $\M_1 \oplus \M_2$ with
the obvious module category structure. A module category is
{\em indecomposable} if it is not equivalent to a
direct sum of two non-trivial module categories.

\begin{comment}
\begin{example}{\bf Indecomposable module categories over pointed categories}.
Let $G$ be a finite group and $\omega \in Z^3(G, \, k^\times)$ be a normalized
$3$-cocycle.
Indecomposable right module categories over $\Vec_G^{\omega}$ correspond to pairs
$(H, \, \psi)$ where $H$ is a subgroup of $G$ such that
$\omega|_{H \times H \times H}$ is cohomologically trivial
and $\psi$ is a $2$-cochain in $C^2(H, \, k^\times)$ satisfying
$\omega|_{H \times H \times H} = d\psi$ \cite{O2}.
Let $\M := \M(H, \, \psi)$ denote the right module category
constructed from the pair $(H, \, \psi)$. The simple objects of
$\M$ are given by the set $H \backslash G$ of right cosets of $H$ in $G$,
the action of $\Vec_G^{\omega}$
on $\M$ comes from the action of $G$ on $H \backslash G$, and the
module category structure isomorphisms are induced from $\psi$.
\end{example}
\end{comment}

Let $\M_1$ and $\M_2$ be two right module categories over
a tensor category $\C$. Let $(F^1, \, \gamma^1)$ and $(F^2, \, \gamma^2)$
be module functors from $\M_1$ to $\M_2$. A {\em natural module
transformation} from $(F^1, \, \gamma^1)$ to $(F^2, \, \gamma^2)$ is a
natural transformation $\eta: F^1 \to F^2$ such that the following
equation holds for all $M \in \M_1$, $X \in \C$:
\begin{equation*}
\label{module trans square}
(\eta_M \otimes \id_X) \, \circ \, \gamma_{M, \, X}^1
= \gamma_{M, \, X}^2 \, \circ \, \eta_{M \otimes X}.
\end{equation*}

Let $\C$ be a tensor category and let $\M$ be a
right module category over $\C$.
The {\em dual category} of $\C$ with respect to
$\M$ is  the category $\C^*_\M:=Fun_\C(\M,\M)$ whose objects are $\C$-module
functors from $\M$ to itself and morphisms are natural
module transformations.
The category $\C^*_\M$ is a tensor category with tensor product being composition
of module functors.
It is known that if $\C$ is a fusion category
and $\M$ is a semisimple $k$-linear indecomposable module category
over $\C$, then $\C^*_\M$ is a fusion category \cite{ENO}.

Two fusion categories $\C$ and $\D$ are said to be {\em weakly Morita
equivalent} if there exists an indecomposable (semisimple $k$-linear) right
module category $\M$ over $\C$ such that the categories $\C^*_{\M}$ and
$\D$ are equivalent as fusion categories.
It was shown by M\"uger \cite{Mu} that this is indeed an equivalence relation.

Consider the fusion category $\Vec_G^{\omega}$, where $G$ is a
finite group and $\omega$ is a normalized $3$-cocycle on $G$. Let
$H$ be a subgroup of $G$ such that $\omega|_{H \times H \times H}$
is cohomologically trivial. Let $\psi$ be a $2$-cochain in $C^2(H,
\, k^\times)$ satisfying $\omega|_{H \times H \times H} = d\psi$.
The twisted group algebra $k^{\psi}[H]$ is an associative unital
algebra in $\Vec_G^{\omega}$. Define $\C = \C(G, \, \omega, \, H,
\, \psi)$ to be the category of $k^{\psi}[H]$-bimodules in
$\Vec_G^{\omega}$. Then $\C$ is a fusion category with tensor
product $\ot_{k^{\psi}[H]}$ and unit object $k^{\psi}[H]$.

Categories of the form $\C(G, \, \omega, \, H, \, \psi)$
are known as {\em group-theoretical} \cite[Definition 8.40]{ENO}, \cite{O2}.
It is known that a fusion category
$\C$ is group-theoretical if and only if it is weakly Morita equivalent
to a pointed category with respect to some indecomposable module category
\cite[Proposition 8.42]{ENO}. More precisely, $\C(G, \, \omega, \, H, \, \psi)$
is equivalent to $(\Vec_G^{\omega})^*_{(H,\psi)}$.

\end{subsection}
%%%%%%%%%%%%%%%%%%%%%%%%%%%%%%%%%%%%%%%
\begin{subsection}{Nilpotent based rings and nilpotent fusion categories}\hfill

Let $\mathbb{Z}_+$ be the semi-ring of non-negative integers.
Let $R$ be a ring
with identity which is a finite rank $\mathbb{Z}$-module. A
{\em $\mathbb{Z}_+$-basis} of $R$ is a basis $B$ such that for all
$X, Y \in B$, $XY = \sum_{Z \in B} n_{X, \, Y}^Z \, Z,$ where
$n_{X, \, Y}^Z \in \mathbb{Z}_+$. An element of $B$ will be called {\em
basic}.

Define a non-degenerate symmetric $\mathbb{Z}$-valued inner
product on $R$ as follows. For all elements $X= \sum_{Z \in B} \, a_Z Z$
and $Y= \sum_{Z \in B} \, b_Z Z$ of $R$ we set
\begin{equation}
\label{pairing}
(X, Y) =\sum_{Z \in B}\, a_Zb_Z.
\end{equation}

\newpage

\begin{definition}[\cite{O1}]
A {\em based ring}
is a pair $(R,\,B)$ consisting of a ring $R$ (with identity $1$) with a
$\mathbb{Z}_+$-basis $B$ satisfying the following properties:
\begin{enumerate}
\item[(1)] $1 \in B$.
\item[(2)] There is an involution $X \mapsto X^*$ of $B$ such that
the induced map \linebreak $X = \sum_{W \in B} a_W W \mapsto X^* =
\sum_{W \in B} a_W W^*$ satisfies
$$(XY,\,Z) = (X, ZY^*) = (Y,\, X^*Z)$$ for all $X,Y,Z\in R$.
\end{enumerate}
\end{definition}

By a {\em based subring} of a based ring $(R,\,B)$ we will mean a
based ring $(S, C)$ where $C$ is a subset  of $B$ and $S$ is a subring of $R$.

Let us recall some definitions from \cite{GN}.

Let $R = (R, \, B)$ be a based ring and let $\C$ be a fusion category.

Let $R_{ad}$ denote the based subring of $R$ generated by all basic
elements of $R$ contained in $XX^*$, $X \in B$. Let $R^{(0)} := R$,
$R^{(1)} := R_{ad}$, and $R^{(i)} := (R^{(i-1)})_{ad}$, for every positive
integer $i$.
Similarly, let $\C_{ad}$ denote the full fusion subcategory of $\C$
generated by all simple subobjects of $X \ot X^*$, $X$ a simple
object of $\C$. Let $\C^{(0)} := \C$,
$\C^{(1)} := \C_{ad}$, and $\C^{(i)} := (\C^{(i-1)})_{ad}$,
for every positive integer $i$.

$R$ is said to be {\em nilpotent} if $R^{(n)} = \mathbb{Z}1$, for some
$n$. The smallest $n$ for which this happens is called the
{\em nilpotency class} of $R$ and is denoted by $cl(R)$.

$\C$ is said to be {\em nilpotent} if $\C^{(n)} \cong \Vec$, for some
$n$. The smallest $n$ for which this happens is called the
{\em nilpotency class} of $\C$ and is denoted by $cl(\C)$.

Note that a fusion category is nilpotent if and only if its
Grothendieck ring is nilpotent. Also note that for any finite group
$G$, the fusion category $\Rep(G)$ of representations of $G$ is
nilpotent if and only if the group $G$ is nilpotent.

Let $\C$ be a fusion category. We can view $\C$ as a $\C_{ad}$-bimodule category.
As such, it decomposes into a direct sum of indecomposable $\C_{ad}$-bimodule categories:
$\C = \oplus_{a \in A} \C_a$, where A is the index set.
It was shown in \cite{GN} that there is a canonical group structure on the index set A.
This group is called the {\em universal grading group} of $\C$ and is denoted by $U(C)$.
Every fusion category is faithfully graded (in the sense of \cite[Definition 5.9]{ENO})
by its universal grading group.

\end{subsection}
%%%%%%%%%%%%%%%%%%%%%%%%%%%%%%%%%%%%%%%%%%%%%%%%%
\begin{subsection}{Some definitions and results from group theory}\hfill
\label{group theory}

The following definitions and results are contained in \cite{R}.

Let $H$ be a subgroup of a group $G$. The subgroup $H$ is said to
be {\em subnormal} in $G$ if there exist subgroups $H_1, \cdots
H_{n-1}$ of $G$ such that
\begin{equation*}
H = H_0 \unlhd H_1 \unlhd \cdots \unlhd H_{n-1} \unlhd H_n = G.
\end{equation*}

For any non-empty subsets $X$ and $Y$ of $G$, let $X^Y$ denote the
subgroup generated by the set $\{yxy^{-1} \mid x \in X, y \in Y \}$.
Define a sequence of subgroups $H^{(G, \,i)}, i = 0, 1, \dots$, of $G$ by the
rules
\begin{equation*}
H^{(G, \, 0)} := G \mbox{  and  } H^{(G, \, i+1)} := H^{H^{(G, \,i)}}.
\end{equation*}
So we get the following sequence
\begin{equation*}
G = H^{(G, \, 0)} \unrhd H^{(G, \, 1)}
\unrhd H^{(G, \, 2)} \unrhd \cdots.
\end{equation*}
Note that $H^{(G, \, 1)}$ is the normal closure of $H$ in $G$. The above
sequence is called the {\em series of successive normal closure}
of $H$ in $G$. It is known that $H$ is subnormal in $G$ if and
only if $H^{(G, \,n)} = H$ for some $n \geq 0$. If $H$ is subnormal in $G$,
the smallest $n$ for which $H^{(G, \, n)} = H$ is called the {\em defect}
of $H$ in $G$.

Suppose $G$ is finite. Then it is known that $G$ is nilpotent if
and only if any subgroup of $G$ is subnormal in $G$. It is also
known that if $H$ is nilpotent and is subnormal in $G$, then the
normal closure of $H$ in $G$ is nilpotent. Indeed, it can be shown
that if $H$ is nilpotent and is subnormal in $G$, then $H$ is
contained in the Fitting subgroup $\mbox{Fit}(G)$ of $G$ (= the
unique largest normal nilpotent subgroup of $G$), and hence the
normal closure of $H$ in $G$ must be nilpotent.
\end{subsection}
%%%%%%%%%%%%%%%%%%%%%%%%%%%%%%%%%%%%%%%%%%%%%%%%%

\end{section}
%%%%%%%%%%%%%%%%%%%%%%%%%%%%%%%%%%%%%%%%%%%%%%%%%%%%%%%%%%%%%%%%%
\begin{section}
{Fusion categories graded by based rings and double coset rings}

In this section we define the notion of a fusion category graded by a based ring
(generalizing the notion of a fusion category graded by a finite group).
We then define the double coset based ring and give a necessary and sufficient condition
for it to be nilpotent.

%%%%%%%%%%%%%%%%%%%%%%%%%%%%%%%%%%%%%%%%%%%%%%%%%
\begin{subsection}{Fusion categories graded by based rings}\hfill

\begin{definition}
A fusion category $\C$ is said to be {\em graded} by a based ring $(R,\, B)$ if
$\C$ decomposes into a direct sum of full abelian subcategories
$\C =\oplus_{X \in B}\,\C_X$ such that
$(\C_X)^* = \C_{X^*}$ and $\C_{X} \ot \C_{Y} \subseteq
\oplus_{Z \in \{W \in B \mid W \text{ is contained in } XY\}} \, \C_{Z}$, for all $X,Y \in B$.
\end{definition}

\begin{remark}
Note that the trivial component $\C_1$ is a fusion subcategory of $\C$.
\end{remark}

Let $\C$ be a fusion category which is graded by a based ring
$(R,\, B)$.

\begin{definition}
For any subcategory $\D \subseteq \C$, define its {\em support}
$\mbox{Supp}(\D) := \{X \in B \mid \D \cap \C_X \not = \{0\} \}$.
We will say that $\C$ is {\em faithfully} graded by $(R,\, B)$
if $\C_X \not = \{0\}$ and
$\mbox{Supp}(\C_{X} \ot \C_{Y}) = \{W \in B \mid W \text{ is contained in } XY\}$,
for all $X, Y \in B$.
\end{definition}

\begin{remark}
(i) Every fusion category is faithfully graded by its Grothendieck ring.\\
(ii) Every fusion category that is graded by a group $G$ is graded by the based ring
$(\mathbb{Z}G, \, G)$.
\end{remark}

Recall that for any fusion category $\C$, $\C_{ad}$ denotes the full
fusion subcategory of $\C$ generated by all simple subobjects of $X \ot X^*$, $X$
a simple object of $\C$; $\C^{(0)} = \C$,
$\C^{(1)} = \C_{ad}$, and $\C^{(i)} = (\C^{(i-1)})_{ad}$
for every positive integer $i$.

Also recall that for any based ring $(R, \, B)$, $R_{ad}$ denotes the based subring
of $R$ generated by all basic elements of $R$ contained in $XX^*$, $X \in B$;
$R^{(0)} = R$, $R^{(1)} = R_{ad}$, and $R^{(i)} = (R^{(i-1)})_{ad}$ for every positive
integer $i$.

\begin{proposition}
\label{nilpotent graded fusion}
Let $\C$ be a fusion category that is faithfully graded by a
based ring $R = (R, \, B)$. Then $\C$ is nilpotent if and only if $R$
is nilpotent and the trivial component $\C_1$ is nilpotent.
If $\C$ is nilpotent, then its nilpotency class $cl(\C)$ satisfies the following
inequality:
$$cl(R) \leq cl(\C) \leq cl(R) + cl(\C_1).$$
\end{proposition}
\begin{proof}
Since the grading of $\C$ by $R$ is faithful, we have
$\mbox{Supp}(\C^{(i)}) = B \cap R^{(i)}$ for any non-negative
integer $i$. Indeed, note that even without faithfulness of the
grading we have $\mbox{Supp}(\C^{(i)}) \subseteq B \cap R^{(i)}$.
Faithfulness of the grading implies that $B \cap R^{(i)} \subseteq
\mbox{Supp}(\C^{(i)})$. Now suppose that $\C$ is nilpotent of
nilpotency class $n$. Then the trivial component $\C_1$ being a
fusion subcategory of $\C$ is nilpotent. Also,
$\mbox{Supp}(\C^{(n)})$ must be equal to $\{1\}$. It follows that
$R$ must be nilpotent. Conversely, suppose that the trivial
component $\C_1$ is nilpotent and $R$ is nilpotent of nilpotency
class $n$. Then $\C^{(n)} \subseteq \C_1$ and it follows that $\C$
must be nilpotent. The statement about nilpotency class should be
evident and the proposition is proved.
\end{proof}

\end{subsection}
%%%%%%%%%%%%%%%%%%%%%%%%%%%%%%%%%%%%%%%%%%%%%%%%%%
\begin{subsection}{The double coset ring}\hfill

Let $H$ be a subgroup of a finite group $G$.
Let $\R(G, \, H)$ denote the free $\mathbb{Z}$-module generated by the set
$\O$ of double cosets of $H$ in $G$.
For any $HxH, HyH \in \O$, the set $HxHyH$ is a union of double cosets.
Define the product $HxH \cdot HyH$ by
$$HxH \cdot HyH := \sum_{HzH \in \O}N_{HxH, \, HyH}^{HzH} \, HzH,$$
where
$$N_{HxH, \, HyH}^{HzH} =
\begin{cases} 1 \text{ if } HzH \subseteq HxHyH,\\
0 \text{ otherwise}.
\end{cases}$$

This multiplication rule on $\O$ extends, by linearity, to a
multiplication rule on $\R(G, \, H)$. The identity element of
$\R(G, \, H)$ is given by the trivial double coset $H = H1_GH$.
There is an involution $*$ on the set $\O$ defined as follows. For
any $HxH \in \O$, define $(HxH)^* := Hx^{-1}H$. It is
straightforward to check that $\R(G,\, H)$ is a based ring.

Let $\S$ be a based subring of $\R(G, \, H)$.
Define
\begin{equation*}
\Gamma_{\S} := \bigcup_{X \in {\S} \cap \O} X.
\end{equation*}
Note that $\Gamma_{\S}$ is a subgroup of $G$ that contains $H$.
Also note that $\Gamma_{\R(G, \, H)} = G$.

\begin{lemma}
The assignment $\S \mapsto \Gamma_{\S}$ is a bijection between the set of
based subrings of the double coset ring $\R(G, H)$ and the set of subgroups
of $G$ containing $H$.
\end{lemma}
\begin{proof}
Let $K$ be a subgroup of $G$ that contains $H$. The double coset ring
$\R(K,\, H)$ is a based subring of $\R(G, \, H)$. It is evident
that the assignment $K \mapsto \R(K, \, H)$ is inverse to the
assignment defined in the statement of the lemma.
\end{proof}

\begin{proposition}
\label{nilpotent double coset}
The double coset  ring $\R(G, H)$ is nilpotent if and only if $H$ is subnormal
in $G$. If $\R(G, H)$ is nilpotent, then its nilpotency class is
equal to the defect of $H$ in $G$.
\end{proposition}
\begin{proof}
Let $\R = \R(G, H)$. Observe that $\Gamma_{\R^{(i)}} = H^{(G, \,
i)}$, for all non-negative integers $i$ (see Subsection \ref{group
theory} for the definition of $H^{(G, \, i)}$). Note that $\R$ is
nilpotent if and only if $H^{(G, \, n)} = H$ for some non-negative
integer $n$. The latter condition is equivalent to the condition
that $H$ is subnormal in $G$. Recall that if $H$ is subnormal in
$G$, then the defect of $H$ in $G$ is defined to be the smallest
non-negative integer $n$ such that $H^{(G, \, n)} = H$. It follows
that if $\R$ is nilpotent, then its nilpotency class is equal to
the defect of $H$ in $G$.
\end{proof}

\end{subsection}

\end{section}
%%%%%%%%%%%%%%%%%%%%%%%%%%%%%%%%%%%%%%%%%%%%%%%%%%%%%%%%%%%%%%%%%%%%%%%%%%%%%
\begin{section}{Nilpotency of a group-theoretical category}\hfill

In this section we give a necessary and sufficient condition for a group-theoretical
category to be nilpotent.

We start with the following theorem.

\begin{theorem}
\label{grading}
Let $\C = \C(G, \, \omega, \, H, \, \psi)$ be a group-theoretical category.
Then $\C$ is faithfully graded by the double coset ring $\R(G, \, H)$, with the trivial
component being the representation category $\Rep(H)$ of $H$.
\end{theorem}
\begin{proof}
It follows from the results in \cite{O2} that the set of isomorphism classes of
simple objects in $\C$ are
parametrized by pairs $(a,\rho)$, where $a\in G$ is a representative of a double
coset $X:=HaH$ of $H$ in $G$ (i.e., a basic element $X$ in $\R(G, \, H)$) and an irreducible
projective representation of $H^a:=H\cap aHa^{-1}$ with a certain $2-$cocycle. Moreover,
the tensor product of two simple objects $X$, $Y$, corresponding to $(a,\rho)$, $(b,\tau)$,
respectively, is supported
on the union of the double cosets appearing in the decomposition of
$XY$. Therefore if we let $\C_X$, $X:=HaH$, be the subcategory of $\C$
generated by all simple objects which correspond to pairs $(a,\rho)$, we get that
$\C=\oplus_{X}\,\C_X$, as required. It is clear that $\C_H=\Rep(H)$.
\end{proof}

\begin{remark}
We note that if $N$ is the normal closure of $H$ in $G$ then the group ring $\mathbb{Z}[G/N]$
is a homomorphic image of $\R(G, \, H)$. Hence the group-theoretical category
$\C = \C(G, \, \omega, \, H, \, \psi)$ is $G/N-$graded.
\end{remark}

\begin{corollary}
\label{main corollary}
Let $\C = \C(G, \, \omega, \, H, \, \psi)$ be a group-theoretical category.
Then $\C$ is nilpotent if and only if the normal closure of $H$ in $G$ is nilpotent.
If $\C$ is nilpotent, then its nilpotency class $cl(\C)$ satisfies the following
inequality:
$$cl(H) \leq cl(\C) \leq cl(H) + (\text{defect of } H \text{ in } G).$$
\end{corollary}
\begin{proof}
By Theorem \ref{grading} and Proposition \ref{nilpotent graded fusion},
it follows that $\C$ is nilpotent if and only if
the double coset ring $\R(G, \, H)$ is nilpotent and $H$ is nilpotent. By
Proposition \ref{nilpotent double coset}, $\R(G, \, H)$ is nilpotent
if and only if $H$ is subnormal in $G$.
Since $G$ is a finite group, it follows from the remarks in Subsection
\ref{group theory} that $H$ is nilpotent and is subnormal
in $G$ if and only if the normal closure of $H$ in $G$ is
nilpotent. The statement about the nilpotency class of $\C$ follows
immediately from Proposition \ref{nilpotent graded fusion} and
Proposition \ref{nilpotent double coset}.
\end{proof}

\begin{example}
Let $G$ be a finite group and let $\omega$ be a $3$-cocycle on
$G$. It was shown in \cite{O2} that the representation category
$\Rep(D^{\omega}(G))$ of the twisted quantum double of $G$ is
equivalent to $\C(G \times G, \, \tilde{\omega}, \, \Delta(G), \,
1)$, where $\tilde{\omega}$ is a certain $3$-cocycle on $G \times
G$ and $\Delta(G)$ is the diagonal subgroup of $G$. It follows
from Corollary \ref{main corollary} that $\Rep(D^{\omega}(G))$ is
nilpotent if and only if $G$ is nilpotent.
\end{example}

\end{section}
%%%%%%%%%%%%%%%%%%%%%%%%%%%%%%%%%%%%%%%%%%%%%%%%%%
\begin{section}{The pointed subcategory of a group-theoretical category}

In this section we describe the simple objects in a
group-theoretical category and then describe the group of
invertible objects in a group-theoretical category.

\begin{subsection}{Simple objects in a group-theoretical category}\hfill

Let $\C = \C(G, \, \omega, \, H, \, \psi)$ be a group-theoretical
category. Let $R = \{u(X) \mid X \in H \backslash G / H\}$ be a
set of representatives of double cosets of $H$ in $G$. We assume
that $u(H1_GH) = 1_G$. In \cite{O2} it is explained how a simple
object in $\C$ gives rise to a pair $(g, \, \overline{\rho})$,
where $g \in R$ and $\overline{\rho}$ is the isomorphism class of
an irreducible projective representation $\rho$ of $H^g$ with a
certain $2$-cocycle $\psi^g$. Let us recall this in details.

For each $g \in G$, let $H^g := H \cap gHg^{-1}$. The group $H^g$
has a well-defined $2$-cocycle $\psi^g$ defined by
\begin{equation*}
\psi^g(h_1, \, h_2) := \psi(h_1, \, h_2) \psi(g^{-1}h_2^{-1}g, \, g^{-1}h_1^{-1}g)
\frac{\omega(h_1, \, h_2, \, g) \omega(h_1, \, h_2g, \, g^{-1}h_2^{-1}g)}
{\omega(h_1h_2g, \, g^{-1}h_2^{-1}g, \, g^{-1}h_1^{-1}g)}.
\end{equation*}

Let $B = \oplus_{g \in G} B_g$ be an object in $\C$. So $B$ is equipped
with isomorphisms $l_{h, \, g} : B_g \xrightarrow{\sim} B_{hg}$
and $r_{g, \, h} : B_g \xrightarrow{\sim} B_{gh}$, $g \in G, h \in
H$. These isomorphisms satisfy the following identities:
$$\omega(h_1, \, h_2, \, g) \psi(h_1, \, h_2) l_{h_1h_2, \, g}
= l_{h_1, \, h_2g} \circ l_{h_2, \, g},$$
$$\psi(h_1, \, h_2) r_{g, \, h_1h_2} =
\omega(g, \, h_1, \, h_2) r_{gh_1, \, h_2} \circ r_{g, \, h_1}$$
and
$$l_{h_1, \, gh_2} \circ r_{g, \, h_2} = \omega(h_1, \, g, \, h_2)
r_{h_1g, \, h_2} \circ l_{h_1, \, g}.$$

The above three identities say that $B$ is a left $k^{\psi}[H]$-module,
$B$ is a right $k^{\psi}[H]$-module, and that the left and right
module structures on $B$ commute, respectively.
It is clear that $B$ is a direct sum of subbimodules
supported on individual double cosets of $H$ in $G$. Suppose $B$ contains
a subbimodule that is supported on a double coset represented by $g$. Then
one get a projective representation $\rho : H^g \to GL(V)$ with $2$-cocycle
$\psi^g$ defined as follows. Let $V := B_g$ and
\begin{equation}
\label{proj rep}
\rho(h) := r_{hg, \, g^{-1}h^{-1}g}
\circ l_{h, \, g}, \qquad h \in H^g.
\end{equation}

The following theorem, stated in \cite{O2}, asserts that the above correspondence
gives a bijection between isomorphism classes of simple objects in
$\C$ and isomorphism classes of pairs $(g, \, \rho)$. We shall give an alternative
proof of the inverse correspondence by a direct computation.

\begin{theorem}\label{simpobj}
The above correspondence defines a bijection between isomorphism
classes of simple objects in $\C$ and isomorphism classes of pairs
$(g, \, \rho)$, where $g\in R$ and $\rho$ is an irreducible projective
representation of $H^g$ with $2-$cocycle $\psi^g$.
\end{theorem}

\begin{proof}
Given a pair $(g, \, \rho)$, where $g \in R$ and $\rho : H^g
\to GL(V)$ is an irreducible projective representation with $2$-cocycle
$\psi^g$, we assign an object $B$ in $\C$ as follows. Let $T$ be a
set of representatives of $H/H^g$. We assume that $1 \in T$. Let
$B := \oplus_{t \in T, k \in H} B_{tgk}$, where each component is
equal to $V$ as a vector space. The right and left module
structures $r$ and $l$, respectively, on $B$ are defined as
follows.
\begin{equation}
\label{right module}
r_{tgk, \, h} : B_{tgk} \xrightarrow{\sim} B_{tgkh}, v \mapsto
\psi(k, \, h) \omega(tg, \, k, \, h)^{-1} v.
\end{equation}
\begin{equation}
\label{left module}
\begin{split}
&l_{h, \, tgk} : B_{tgk} \xrightarrow{\sim} B_{sg(g^{-1}pg)h},
v \mapsto \frac{\psi(h, \, t)}{\psi(s, \, p) \psi(g^{-1}p^{-1}g, \, g^{-1}pgk)}\\
& \times \frac{\omega(h, \, tg, \, k)\omega(s, \, g, \, g^{-1}pg)
\omega(h, \, t, \, g)} {\omega(s, \, p, \, g)} \\
&\times \frac{\omega(g, \, g^{-1}pg, \, g^{-1}p^{-1}g)
\omega(g^{-1}pg, \, g^{-1}p^{-1}g, \, g^{-1}pgk)}
{\omega(sg, \, g^{-1}pg, \, k)}\rho(p)(v),
\end{split}
\end{equation}
where $s \in T$ and $p \in H^g$ are uniquely determined by the
equation $ht = sp$. It is now straightforward to check that $B$ is
simple, and that the two correspondences are inverse to each
other.
\end{proof}

\end{subsection}

\begin{subsection}{The group of invertible objects in a group-theoretical category}\hfill

For any $g \in N_G(H)$ and $f \in C^n(H, \, k^\times)$, define
${}^g{f} \in C^n(H, \, k^\times)$ by $${}^gf(h_1, \cdots, h_n) :=
f(g^{-1}h_1g, \cdots, g^{-1}h_ng).$$ Pick any $g_1, g_2 \in
N_G(H)$ and let $g_3 = g_1g_2k, k \in H$. Define
\begin{equation}
\label{beta}
\begin{split}
&\beta(g_1, \, g_2) : H \to k^\times, h \mapsto
\frac{\psi(g_2^{-1}g_1^{-1}hg_1g_2k, \, g_3^{-1}h^{-1}g_3)}
{\psi(g_1^{-1}h^{-1}g_1, \, g_1^{-1}hg_1)
\psi(g_2^{-1}g_1^{-1}h^{-1}g_1g_2, \, g_2^{-1}g_1^{-1}hg_1g_2k)}\\
&\times \frac{\omega(g_1^{-1}hg_1, \, g_1^{-1}h^{-1}g_1, \,
g_1^{-1}hg_1) \omega(g_1, \, g_1^{-1}hg_1, \, g_1^{-1}h^{-1}g_1)
\omega(g_1^{-1}hg_1, \, g_2, \, k)}
{\omega(g_2, \, g_2^{-1}g_1^{-1}hg_1g_2, \, k)}\\
&\times \frac{\omega(g_2^{-1}g_1^{-1}hg_1g_2, \,
g_2^{-1}g_1^{-1}h^{-1}g_1g_2, g_2^{-1}g_1^{-1}hg_1g_2k)
\omega(g_2, \, g_2^{-1}g_1^{-1}hg_1g_2, \,
g_2^{-1}g_1^{-1}h^{-1}g_1g_2)} {\omega(g_2, \,
g_2^{-1}g_1^{-1}hg_1g_2k, \, g_3^{-1}h^{-1}g_3)}.
\end{split}
\end{equation}

It is straightforward (but tedious) to verify that
\begin{equation}
\label{beta relation} \psi^{g_3} = d(\beta(g_1, \, g_2)) \,
\psi^{g_1} \, ({}^{g_1}(\psi^{g_2})).
\end{equation}

Let $K := \{g \in R \mid g \in N_G(H) \text{ and } \psi^g
\text{ is cohomologically trivial} \}$. For any $g_1, g_2 \in K$,
define $g_1 \cdot g_2 := u(g_1g_2)$. It follows from \eqref{beta relation}
that with this product rule $K$ is a group that is isomorphic to a subgroup
of $N_G(H)/H$.

For each $g \in K$, fix $\eta_g : H \to k^\times$ such that $d\eta_g = \psi^g$.
We take $\eta_1 := \beta(1, \, 1)^{-1}$.
For any $g_1, g_2 \in K$, define
\begin{equation}
\label{nu}
\nu(g_1, \, g_2) := \frac{\eta_{g_1} ({}^{g_1}\eta_{g_2})}{\eta_{g_1 \cdot g_2}}
\beta(g_1, \, g_2).
\end{equation}

Let $\widehat{H} := \Hom(H, \, k^\times)$ and define a group $K
\ltimes_{\nu} \widehat{H}$ as follows. As a set \linebreak $K
\ltimes_{\nu} \widehat{H} = K \times \widehat{H}$ and for any
$(g_1, \, \rho_1), (g_2, \, \rho_2) \in K \ltimes_{\nu}
\widehat{H}$, define $$(g_1, \, \rho_1) \cdot (g_2, \, \rho_2) =
(g_1 \cdot g_2, \, \nu(g_1, \, g_2) \rho_1 ({}^{g_1}\rho_2)).$$

\begin{theorem}\label{grinv}
The group $G(\C)$ of isomorphism classes of invertible objects of $\C$ is isomorphic to
the group $K \ltimes_{\nu} \widehat{H}$ constructed above.
\end{theorem}
\begin{proof}
By Theorem \ref{simpobj}, $G(\C)$ is in bijection with the set $$L = \{(g,
\, \rho) \mid g \in K, \rho:H \to k^\times \text{ such that }
d\rho = \psi^g \}.$$ The set $L$ becomes a group with product
$$(g_1, \, \rho_1) \cdot (g_2, \, \rho_2) = (g_1 \cdot g_2, \,
\beta(g_1, \, g_2) \rho_1 ({}^{g_1}\rho_2)).$$ The identity
element of $L$ is $(1, \, \beta(1, \, 1)^{-1})$. Let $B, B^\prime$
be objects in $\C$ corresponding to $(g_1, \, \rho_1), (g_2, \,
\rho_2) \in L$, respectively. So $B = \oplus_{h \in H} k_{g_1h}$
and $B^\prime = \oplus_{h \in H} k_{g_2h}$, where each component
is equal to the ground field $k$. The right and left module
structures on $B, B^\prime$ are defined via \eqref{right module}
and \eqref{left module}. Let $A := k^{\psi}[H]$. We have $B \ot_A
B^\prime = (k_{g_1} A) \ot_A (\oplus_{h \in H} k_{g_2h}) = k_{g_1}
\ot (\oplus_{h \in H} k_{g_2h})$. Taking into account \eqref{right
module} and \eqref{left module} we calculate that the projective
representation (defined in \eqref{proj rep}) $\rho: H \to
k^\times$ with $2$-cocycle $\psi^{g_3}$, corresponding to $B \ot_A
B^\prime$, where $g_3 = g_1 \cdot g_2$, is given by $\beta(g_1, \,
g_2) \rho_1 ({}^{g_1} \rho_2)$. So $G(\C)$ is isomorphic to the
group $L$. The map $L \to K \ltimes_{\nu} \widehat{H} : (g, \,
\rho) \mapsto (g, \, \eta_g^{-1} \rho)$ establishes the desired
isomorphism and the theorem is proved.
\end{proof}

\end{subsection}

\end{section}
%%%%%%%%%%%%%%%%%%%%%%%%%%%%%%%%%%%%%%%%%%%%%%%%%%
\begin{section}{The universal grading group of certain group-theoretical categories}\hfill

Recall that every fusion category $\C$ is faithfully graded by its
universal grading group $U(\C)$: $\C = \oplus_{x \in U(\C)} \C_x$.
In this section we describe $U(\C)$ for certain group-theoretical categories.

\begin{lemma}
\label{univ injection} Let $\D$ be a fusion category and let $\E$
be a fusion subcategory of $\D$. The map $U(\E) \to U(\D)$ defined
by the rule $x \mapsto y$ if and only if $\E_{x} \subseteq \D_y
\cap \E$ is a homomorphism. This homomorphism is injective if and
only if $\D_{ad} \cap \E = \E_{ad}$.
\end{lemma}
\begin{proof}
We have universal gradings: $\D = \oplus_{y \in U(\D)} \D_y$ and
$\E = \oplus_{x \in U(\E)} \E_x$. From the former grading we obtain
$\E = \D \cap \E = \oplus_{y \in U(\D)} (D_y \cap \E)$. Note that this
grading need not be faithful. Since $\E_{ad} \subseteq \D_{ad} \cap \E$,
each component $D_y \cap \E$ is a $\E_{ad}$-submodule category of $\E$.
So, for every $x \in U(\E)$ there is a unique $y \in U(\D)$ such that
$\E_x \subseteq \D_y$. This gives rise to a homomorphism $U(\E) \to U(\D)$. It is
evident that this homomorphism is injective if and only if
$\D_{ad} \cap \E = \E_{ad}$.
\end{proof}

\begin{lemma}
\label{univ group}
The universal grading group $U(\Rep(K))$ of the representation category
of a finite group $K$ is isomorphic to the center $Z(K)$ of $K$.
\end{lemma}
\begin{proof}
This is a special case of Theorem 3.8 in \cite{GN} ($H$ being
the group algebra of $K$).
\end{proof}

\begin{proposition}\label{unigr}
Let $\C = \C(G, \, 1, \, H, \, 1)$. Suppose $H$ is normal in $G$.
Then there is a split exact sequence $1 \to Z(H) \to U(\C) \to G/H \to 1$.
Therefore, $U(\C)$ is isomorphic to the semi-direct product $G/H \ltimes Z(H)$.
\end{proposition}
\begin{proof}
By Theorem \ref{grading}, we have a grading of $\C$ by the group
$G/H$: $\C = \oplus_{x \in G/H} \C^x$, where $\C^x$ is the full
abelian subcategory of $\C$ consisting of objects supported on the
coset $x$. Let $\E := \C^1$. We will first show that $\C_{ad} =
\E_{ad}$. Let $R$ be a set representatives of cosets of $H$ in
$G$. Recall that simple objects of $\C$ correspond to pairs $(a,
\, \rho)$, where $a \in R$ and $\rho$ is an irreducible
representation of $H$. Let $B$ be the object in $\C$ corresponding
to $(a, \, \rho)$ defined via \eqref{right module} and \eqref{left
module}. The dual object $B^*$ corresponds to the pair $(b, \,
({}^b \rho)^*)$, where $b \in R$ is the representative of the
coset $a^{-1}H$. The representation (defined in \eqref{proj rep})
corresponding to $B \ot_{k[H]} B^*$ is given by $\rho \ot
{}^a(({}^b \rho)^*) \cong \rho \ot \rho^*$. This establishes the
equality $\C_{ad} = \E_{ad}$.

By Theorem \ref{grading}, $\E \cong \Rep(H)$ and Lemma \ref{univ
group} implies that $U(\E) \cong Z(H)$. By Lemma \ref{univ
injection}, we get an injective homomorphism $i: Z(H) \to U(\C)$.
From \cite[Corollary 3.7]{GN} we get a surjective homomorphism $p:
U(\C) \to G/H$ which is defined as follows. Note that $\E$
contains $\C_{ad}$. Therefore, each $\C^x$ is a
$\C_{ad}$-submodule category of $\C$. So, for every $y \in U(\C)$
there is a unique $p(y) \in G/H$ such that the component $\C_y$ of
the universal grading $\C = \oplus_{z \in U(\C)} \C_z$ is
contained in $\C^{p(y)}$.

We claim that the sequence $1 \to Z(H) \xrightarrow{i} U(\C)
\xrightarrow{p} G/H \to 1$ is exact. We have $\C_{ad} = \E_{ad}
\cong \Rep(H)_{ad} \cong \Rep(H/Z(H))$. By \cite[Proposition
8.20]{ENO}, it follows that $|U(\C)| = |Z(H)| \frac{|G|}{|H|}$ and
therefore $|\Ker p| = |Z(H)|$. So, it suffices to show that $\Ker
p \subseteq \Im i$. We have $\Ker p = \{ y \in U(\C) \mid \C_y
\subseteq \E \}$. Pick any $y \in \Ker p$ and let $K := \{y \in
U(\C) \mid \C_y \cap \E \not = \{0\}\}$. Then $\E = \oplus_{k \in
K} (\C_k \cap \E)$ is a faithful grading of $\E$. Note that $y \in
K$. By \cite[Corollary 3.7]{GN}, there exists $z \in U(\E)$ such
that $\E_z \subseteq \C_y$, i.e., $y \in \Im i$. This establishes
the exactness of the aforementioned sequence.

Finally, we show that the aforementioned sequence splits. Let $\D$ be the full
fusion subcategory of $\C$ generated by simple objects in $\C$ corresponding
to pairs $(a, \, \rho_0)$, where $a \in R$ and $\rho_0$ is the trivial representation
of $H$. Note that $\D \cong \Vec_{G/H}$ and $U(\D) \cong G/H$. Also note that
$\C_{ad} \cap \D = D_{ad} \cong \Vec$. So, by Lemma \ref{univ injection} we
obtain an injection $j: G/H \to U(\C)$. We claim that $p \circ j = \id_{G/H}$.
Pick any $x \in G/H$ and let $j(x) = y$, i.e., $\D_x \subseteq \C_y$.
We have $\C_y \subseteq \C^{p(y)}$ which implies that $D_x \subseteq \C^{p(y)}$.
It follows that $p(y) = x$ and the proposition is proved.
\end{proof}

\end{section}
%%%%%%%%%%%%%%%%%%%%%%%%%%%%%%%%%%%%%%%%%%%%%%%%%%%%%%%%%%%%%%%%%%%%%%%%%%%%%

\bibliographystyle{ams-alpha}

\end{document}